\documentclass[
preprint,
showkeys,
 amsmath,amssymb,
 aps,
]{revtex4-2}

\usepackage{graphicx}
\usepackage{dcolumn}
\usepackage{bm}
\usepackage[margin=1in]{geometry}
\usepackage{natbib}
\usepackage{graphicx}
\usepackage{xcolor}
\usepackage{amsmath}
\usepackage{caption}
\usepackage{schemata}
\usepackage{amsfonts}
\usepackage{float}
\usepackage{subcaption}
\usepackage{float}
\makeatletter
\let\newfloat\newfloat@ltx
\makeatother

\usepackage{algorithm}
\usepackage{algpseudocode}
\usepackage{hyperref}

\begin{document}


\title{Recovering discrete delayed fractional equations from trajectories}

\author{J. Alberto Conejero}
\email{aconejero@upv.es}
\affiliation{Instituto Universitario de Matem\'atica Pura y Aplicada\\ 
Universitat Polit\`ecnica de Val\`encia, 46022 Val\`encia, Spain.}

\author{\`Oscar Garibo-i-Orts}%
\altaffiliation[Also at ]{GRID - Grupo de Investigación en Ciencia de Datos\\ Valencian International University - VIU, Carrer Pintor Sorolla 21, 46002 Val\`encia, Spain.}
\affiliation{Instituto Universitario de Matem\'atica Pura y Aplicada\\ 
Universitat Polit\`ecnica de Val\`encia, 46022 Val\`encia, Spain.}

\author{Carlos Lizama}
\affiliation{Departamento de Matemática y Ciencia de la Computación \\
Universidad de Santiago de Chile, Las Sophoras 173, Estaci\'on Central, Santiago, Chile.}

\date{October 2022}

\begin{abstract}
We show how machine learning methods can unveil the fractional and delayed nature of discrete dynamical systems. In particular, we study the case of the fractional delayed logistic map. We show that given a trajectory, we can detect if it has some delay effect or not, and also to characterize the fractional component of the underlying generation model. 
\end{abstract}

\maketitle


\section{Introduction}

Discrete Fractional Calculus (DFC) permits a description of dynamical systems of real-world processes with discrete nature where there exist long-term connections between different time steps. Such connections are introduced as discrete-memory effects in the model, and DFC has been shown a reliable approach in order to capturing them. Recent examples can be found in the study of chaos \cite{wu_baleanu2014discrete}, variable-order discrete-time recurrent neural networks \cite{huang2020variable}, and fuzzy systems on discrete time domains \cite{huang2021discrete}. As a matter of fact, there have been found interesting applications in different settings such as image encryption \cite{bai2018novel,liu2013numerical,wu2016image,wu2019new} or earthquakes intensity modeling \cite{kong2021modeling,cristofaro2022fractional}, and \cite{ye2023global}.

The origins of DFC are dated to 1956, when Kutter \cite{Ku56} mentioned for the first time differences of fractional order. In 1974, Diaz and Osler \cite{DiOs74} introduced a discrete fractional difference operator defined as an infinite series. In 1988, Grey and Zhang \cite{Gr-Zh88} developed a fractional calculus for the discrete $\nabla$ (backward) difference operator. Miller and Ross \cite{MR1989} defined a fractional sum via the solution of a linear difference equation. Their definition is the discrete analog of the Riemann-Liouville fractional integral, which can be obtained via the solution of a linear differential equation. In 2007, Atici and Eloe \cite{AE2007} introduced the Riemann-Liouville like fractional difference by using the definition of a fractional sum of Miller and Ross, and developed some of its properties that allow one to obtain solutions of certain fractional difference equations. The presence of chaos in these models was first studied for the logistic map \cite{wu_baleanu2014discrete}; see also \cite{wu2014discrete} for sine maps. 

Later, Wu and Baleanu also studied the chaos in a delayed version of the logistic map. They started from an expression where the forward Euler operator $\Delta$ was equal to the nonlinear right term of the logistic, that is 
$\Delta x(n):=x(n+1)-x(n)$, and then they replaced the Euler operator with the left Caputo discrete difference operator $\Delta^\alpha$, obtaining 

\begin{equation}\label{eq:wu_baleanu_formula}
x(n) =x(0)+\frac{\mu}{\Gamma(\nu)}\sum_{j=1}^n\frac{\Gamma(n-j+\nu)}{\Gamma(n-j+1)}x(j-1)(1-x(j-1)),
\end{equation}
where $\mu$ is a parameter and $\nu$ is a scaling factor. Such a relationship can also be obtained by convolution, using the Cesàro numbers of order $\nu$, $k^{\nu}(j)=\frac{\Gamma(\nu+j)}{\Gamma(\nu)\Gamma(j+1)}$ with $j\in\mathbb{N}_0$, as a memory kernel  in terms of the scaling factor $\nu$ \cite{conejero2019visibility}. 
It is worth mentioning that new fractional difference equations have been introduced with a general logarithm function on time scales as a kernel function. Such equations are known as Caputo–Hadamard fractional difference equations \cite{wu2022caputo,song2022hadamard}.

After their first work on discrete fractional chaos, Wu and Baleanu later introduced a two-dimensional dynamical system including a delay term in \cite{wu_baleanu2015fractional}, that reads as

\begin{equation}
\begin{array}{ll}\label{eq:wu_baleanu_delayed_formula}
x(n) & =\displaystyle x(0)+\frac{\mu}{\Gamma(\nu)}\sum_{j=1}^n\frac{\Gamma(n-j+\nu)}{\Gamma(n-j+1)}x(j-1)(1-y(j-1))\\
y(n) & =x(n-1).
\end{array}
\end{equation}

\begin{figure}[h]
 \centering
 \begin{subfigure}{.45\textwidth}
 \includegraphics[width=\linewidth]{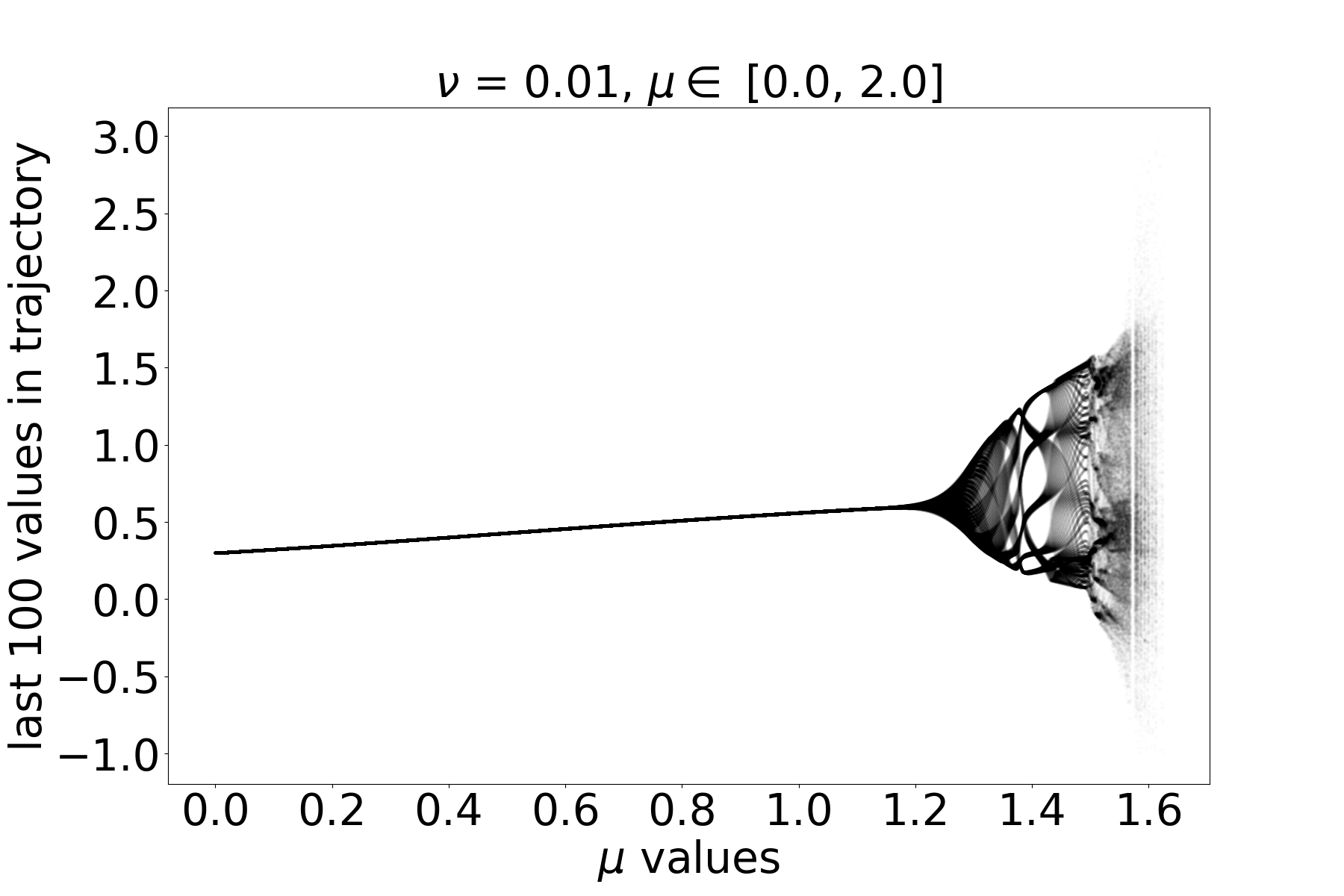}
 \end{subfigure}
 \begin{subfigure}{.45\textwidth}
 \includegraphics[width=\linewidth]{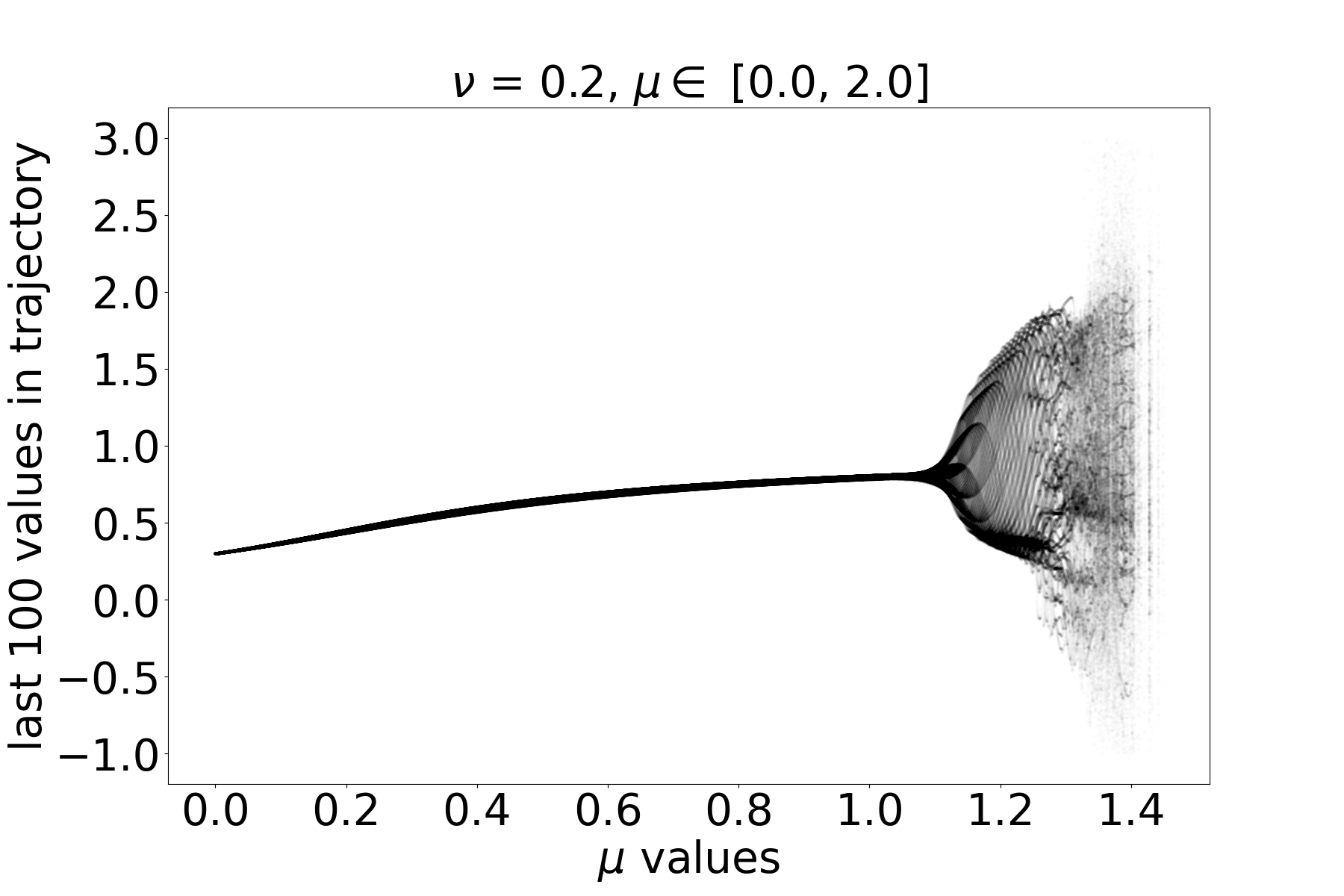}
 \end{subfigure} \\
 \begin{subfigure}{.45\textwidth}
 \includegraphics[width=\linewidth]{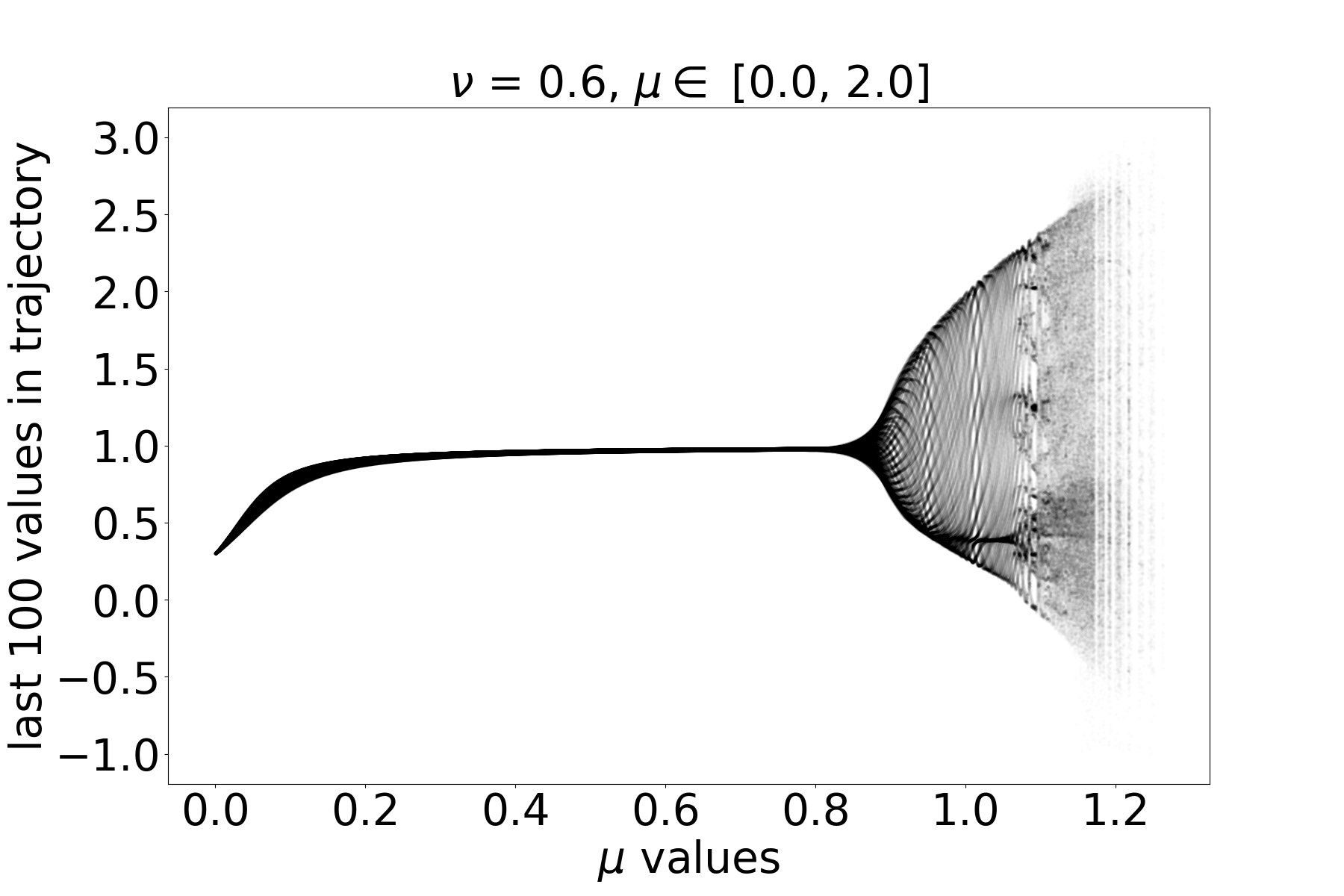}
 \end{subfigure}
 \begin{subfigure}{.45
 \textwidth}
 \includegraphics[width=\linewidth]{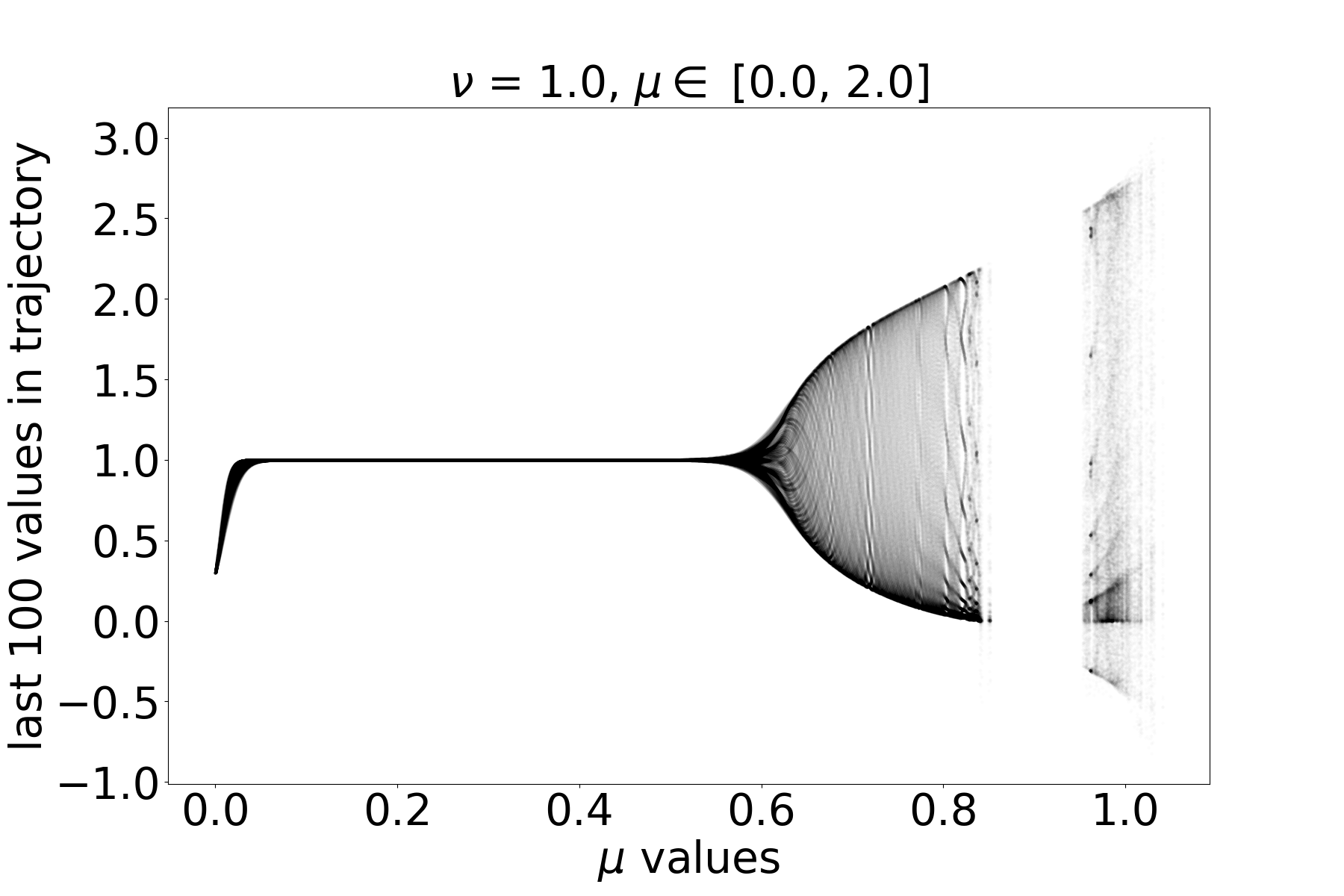}
 \end{subfigure}
 \caption{Feigenbaum plots for the dynamical system given by \eqref{eq:wu_baleanu_formula} for $x(0)=y(0)=0.3$
 and $\nu=0.01$ (top left), $\nu=0.2$ (top right), $\nu=0.6$ (bottom left), and $\nu=1$ (bottom right),
 For each value $\mu$, we compute 200 terms of the sequence, and we plot the last 100 values.}
 \label{fig:feigenbaum_u0_03}
 \end{figure}

\begin{figure}[h]
 \centering
 \begin{subfigure}{.45\textwidth}
 \includegraphics[width=\linewidth]{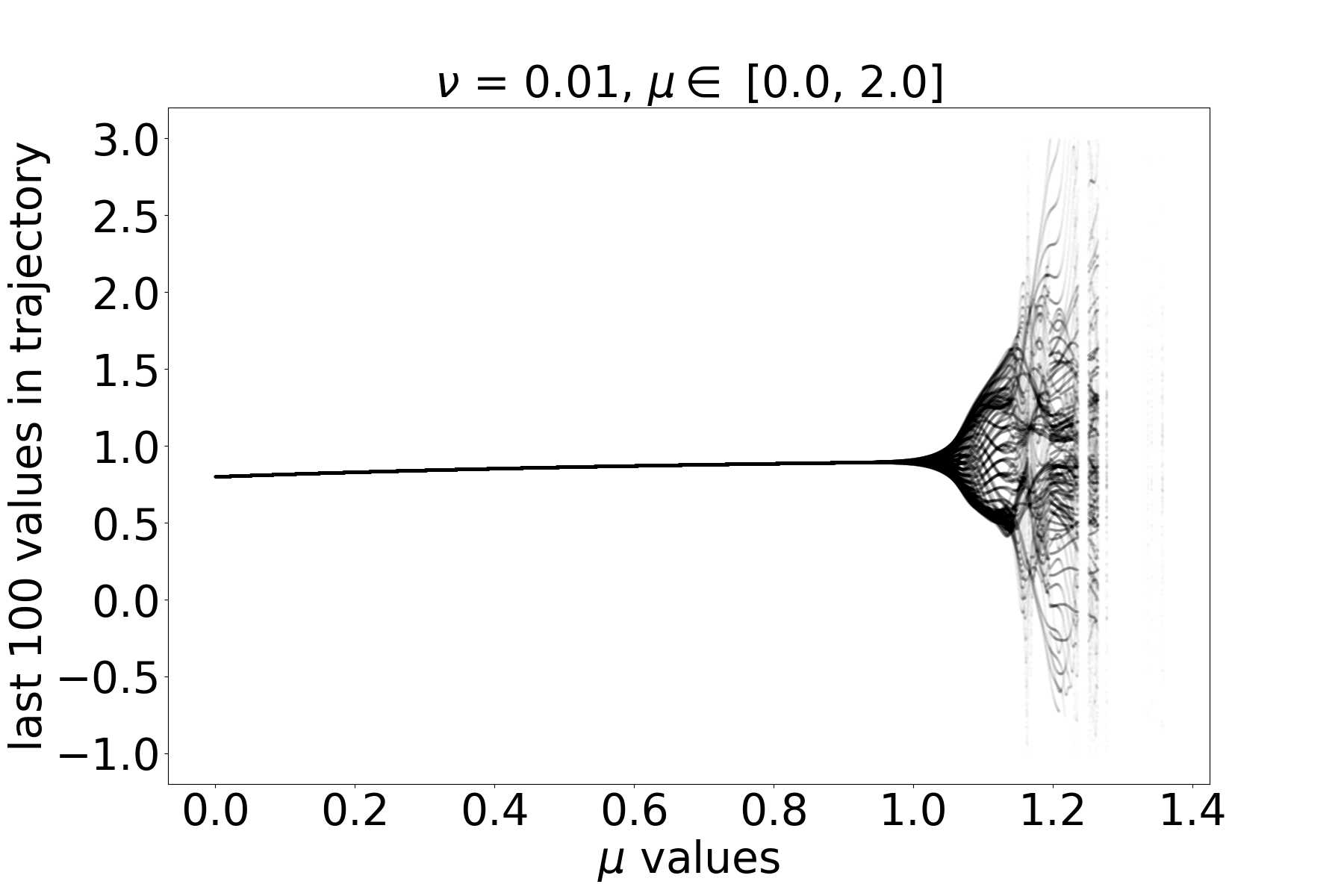}
 \end{subfigure}
 \begin{subfigure}{.45\textwidth}
 \includegraphics[width=\linewidth]{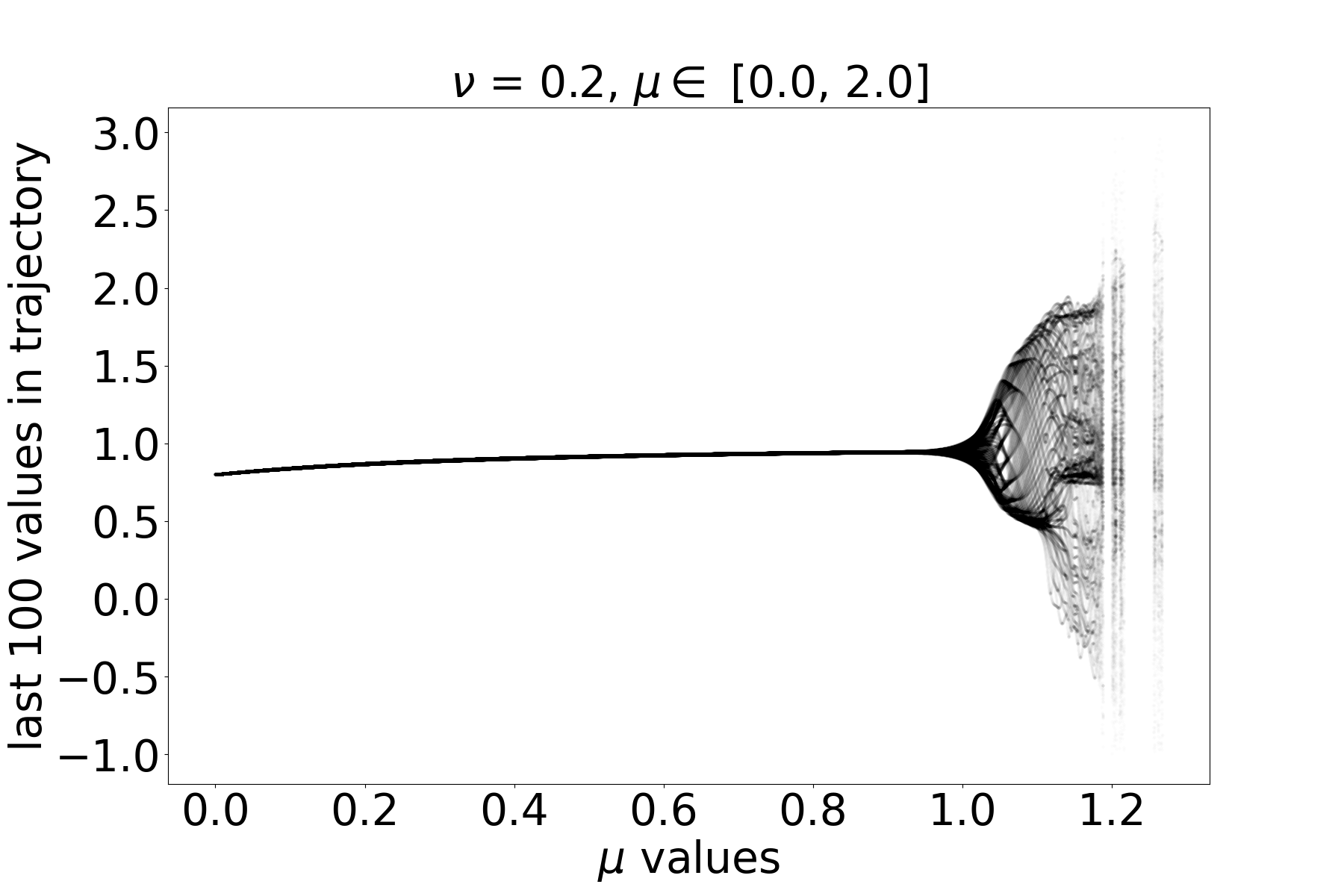}
 \end{subfigure} \\
 \begin{subfigure}{.45\textwidth}
 \includegraphics[width=\linewidth]{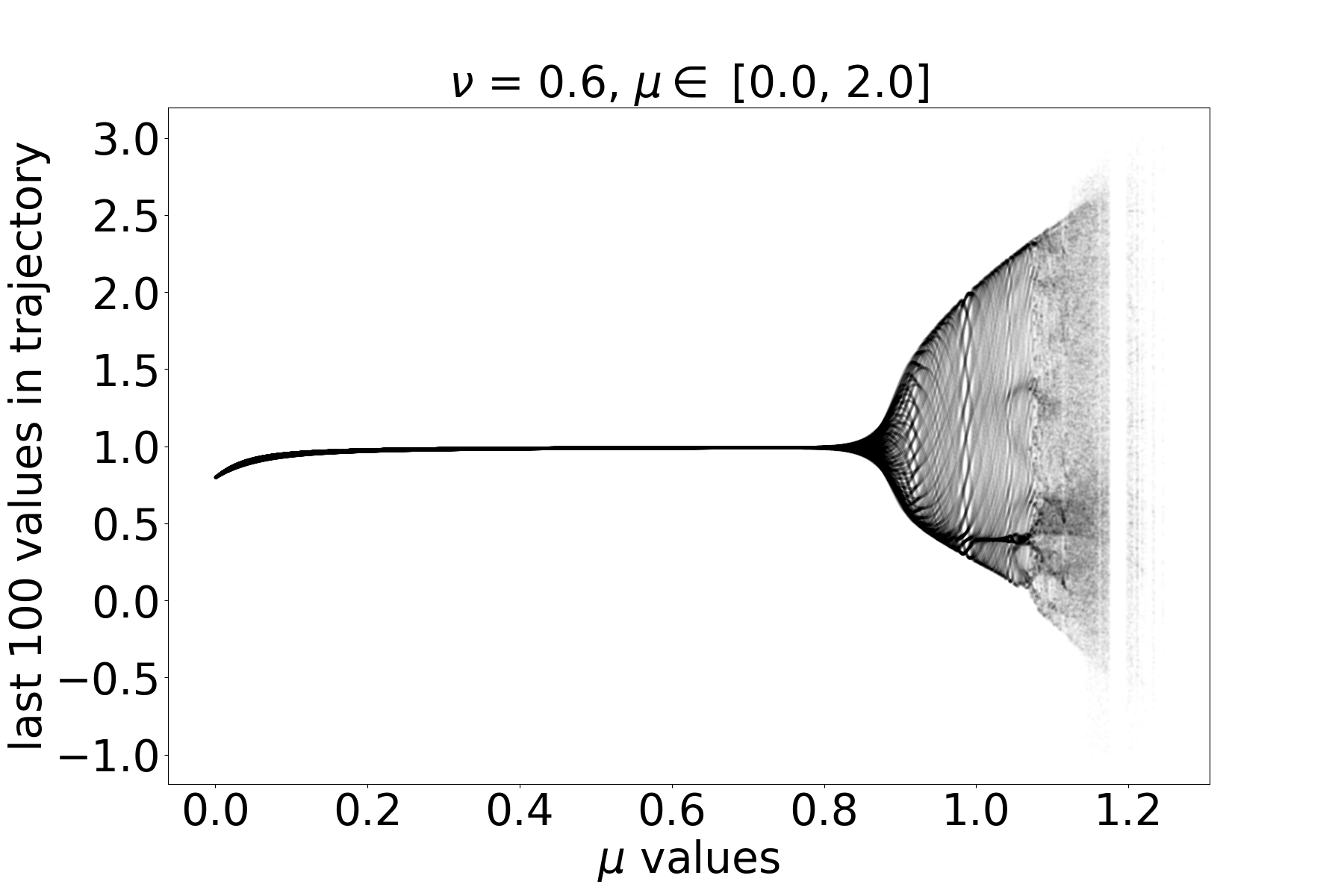}
 \end{subfigure}
 \begin{subfigure}{.45\textwidth}
 \includegraphics[width=\linewidth]{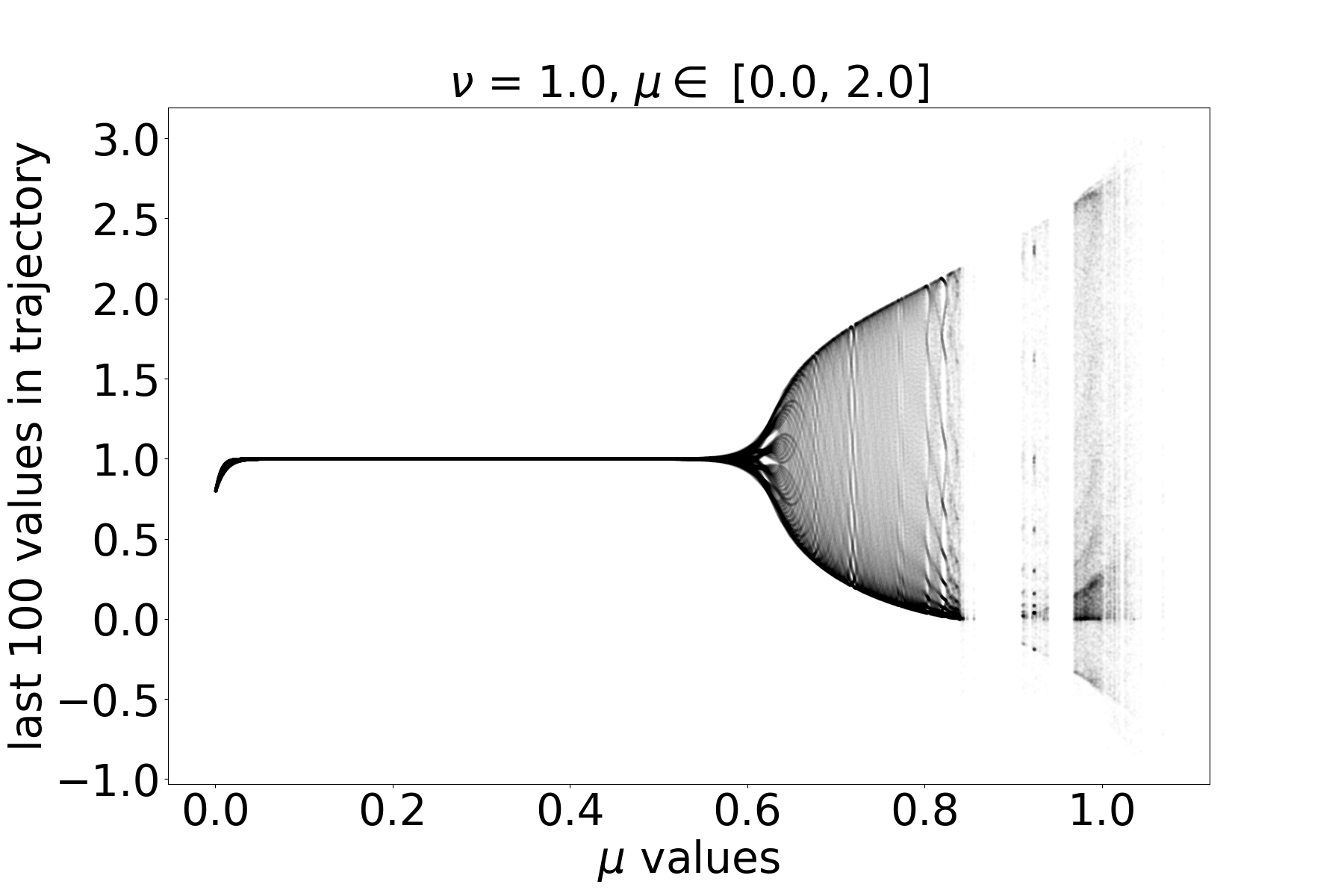}
 \end{subfigure}
 \caption{Feigenbaum plots for the dynamical system given by \eqref{eq:wu_baleanu_formula} for $x(0)=y(0)=0.8$
 and $\nu=0.01$ (top left), $\nu=0.2$ (top right), $\nu=0.6$ (bottom left), and $\nu=1$ (bottom right),
 For each value $\mu$, we compute 200 terms of the sequence, and we plot the last 100 values.}
 \label{fig:feigenbaum_u0_08}
 \end{figure}

This delayed fractional discrete dynamical system also presents chaos as can be noticed in the Feigenbaum diagrams computed in terms of the parameter $\mu$. As a matter of fact, in Figures \ref{fig:feigenbaum_u0_03} \ref{fig:feigenbaum_u0_08}, we show the Feigenbaum diagrams computed for the initial conditions $x(0)=y(0)=0.3$ and $x(0)=y(0)=0.8$. 

The emergence of machine learning as a research tool has been extended to almost all research fields.
It has shown its potential in Physics being used for studying anomalous diffusion noisy trajectories, 
which are characterized as the ones whose variance of the mean square displacement grows with respect to the time in terms of $t^\alpha$, with $\alpha$ being the fractional diffusion exponent \cite{munoz-gil2021etai}. 
Examples of models generating those types of trajectories are Annealed Transient Time Motion (ATTM) \cite{massignan2014nonergodic}, Continuous Time Random Walk (CTRW) \cite{scher1975anomalous}, Fractional Brownian Motion (FBM) \cite{mandelbrot1968fractional,jeon2010fractional}, Lévy Walks (LW) \cite{munoz-gil2021objective, klafter1994levy}, and Scaled Brownian Motion (SBM) \cite{lim2002self-similar}.

Several models have been developed in the frame of the Andi Challenge for inferring the fractional exponent $\alpha$ and for determining the generating model of 1D, 2D, and 3D noisy trajectories, see for instance \cite{argun2021classification,garibo-i-orts2021efficient,gentili2021characterization}. We refer the reader for a full comparison between these models to
\cite{munoz-gil2021objective}.\medskip

Machine learning and artificial intelligence  methods have been also successfully incorporated in the study of formal mathematical problems, as is the case of finding the solution of nonlinear models \cite{raja2018new,sabir2021design} or the recent success of discovering new multiplication algorithms \cite{fawzi2022discovering}.\medskip

Fractional models are suitable for representing dynamical systems with a memory effect, either in space or time, due to the non-local character of fractional operators, as it can be shown in the following examples \cite{baleanu2010new,conejero2022fractional,ilhan2020generalization,sun2018new}. There is not a universal definition of fractional derivative \cite{valerio2022many} and in some cases, it is not straightforward to know which definition would be better to choose \cite{ortigueira2017derivative}.\medskip

Recently, we have studied how machine learning models can be used in connection with fractional models in order to infer the $\mu$ 
parameter and the scaling factor $\nu$ of the fractional version of the logistic equation \cite{conejero_garibo_lizama2022inferring}. In this work, we study up to which point  we can predict the $\mu$ and $\nu$ parameters from short trajectories generated by the discrete fractional delayed dynamical system given by \eqref{eq:wu_baleanu_delayed_formula}. We also analyze if we can discriminate if trajectories contained a delayed component by comparing trajectories generated by the aforementioned systems described in \eqref{eq:wu_baleanu_delayed_formula} and \eqref{eq:wu_baleanu_formula}.

For this purpose, we have used a model that combines convolutional and recurrent neural networks, that presented successful results for inferring the exponent $\alpha$ of anomalous diffusion trajectories \cite{garibo-i-orts2021efficient} and for predicting the number of new COVID19 cases \cite{lozano2021open}. We will show in which cases we can clearly determine from a given trajectory (1) which were the parameters used for generating them, and (2) if we can affirm that it contains a delay component or not.

In Section \ref{sec:methodology}, we briefly outline the architecture of the machine learning model used. We also show how we construct the training, validation, and test data sets. We introduce the results in Section \ref{sec:results}, and we draw some conclusions in Section \ref{sec:conclusions}.\medskip

\section{Methodology}\label{sec:methodology}

\begin{figure}[h!]
	  	\centering	
	  	\includegraphics[width = 1.0\linewidth]{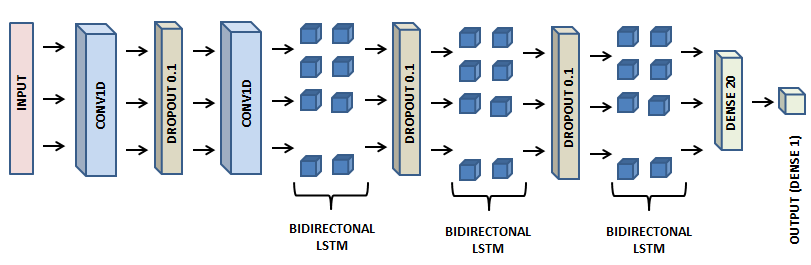}
	  	\caption{Machine learning architecture used for inferring the generating $\mu$ and $\nu$ parameters of the delayed logistic model trajectories.}
	  	\label{fig:architecture}
	 \end{figure}

We start by describing our network architecture we depict in Figure \ref{fig:architecture} and that consists of 3 parts:

\begin{enumerate}
    \item First, a trajectory is processed by two convolutional layers. They are intended for feature extraction since they retain spatial structure inherent to trajectories by setting a patch, which behaves as a sliding window, that goes along the complete trajectory, connecting each patch of size kernel size, to a single neuron. By doing so, different neurons specialize in different regions in the trajectory \cite{krizhevsky2010cifar,krizhevsky2012alexnet}.\medskip
    
    The first convolutional layer consists of 32 filters, and the second one of 64, with a sliding window (kernel) of size 5, following the doubling rule \cite{krizhevsky2010cifar,krizhevsky2012alexnet}.  
    The different number of filters returns features at different scales while retaining the spatial structure.\medskip
    
    \item Secondly, the features feed 3 stacked bidirectional Long Short Term Memory (LSTM) of 64 inputs and 32 units each one \cite{hochreiter1997lstm,lipton2015critical} using Adam optimizer \cite{kingma2014adam}. This type of recurrent neural network is able to learn sequential dependencies from the extracted features.
    For this purpose, LSTM layers include gated cells to avoid short-term memory or gradient vanishing, allowing the layer to learn from the complete trajectory, even from positions at the beginning of it. A cell state channel is selectively updated to allow getting rid of irrelevant information while retaining the important pieces.\medskip
    
    In other words, they are able to retain  memory about early positions in the trajectory, thus considering the complete trajectory while adjusting to the target ($\mu$ and $\nu$). We also include a dropout layer after each of these LSTM of the $10\%$ neurons to avoid over-fitting.\medskip
    \item Lastly, the output of the LSTM blocks is passed to two dense layers of 20 and 1 unit. This last choice depends if we want to predict a single parameter or both of them at the same time. The activation function of this last dense layer depends on the problem to be addressed: It is a linear  function for regression of the  $\mu$ and $\nu$ parameters and a sigmoid function for determining whether a trajectory comes from a delayed model or not, that is, for the classification problem.\medskip
\end{enumerate}

In order to train this model, we built two separate data sets for training and validating while training, the previous model. Once the model was trained, we also constructed an additional data set for testing. These data sets are generated according to the following indications:

\begin{itemize}
    \item $\mu \in[0.0,2.0]$ with increments of $0.001$.
    \item $\nu \in [0.01,1]$ with increments of $0.01$.
    \item trajectory length $N$, with $N\in [10,50]$ randomly selected.
    \item $u_0 \in[0,1]$ randomly chosen with a resolution of $10^{-2}$.
\end{itemize}

We point out that in order to capture the chaotic dynamics that appear in some regions of the $\mu$ parameter (see Figures \ref{fig:feigenbaum_u0_03} and \ref{fig:feigenbaum_u0_08}), the resolution of $\mu$ is increased by 10 times respect to the one of $\nu$. We summarize the data generation procedure in the algorithm \ref{alg:datasetcreation}.\medskip

\begin{algorithm}
  \caption{Data sets creation algorithm}\label{alg:datasetcreation}
  \begin{algorithmic}[1]
    \For{\texttt{<every $\mu$>}}
        \For{\texttt{<every $\nu$>}}
            \State \texttt{<i = 0>}
            \While{i $<$ 5}
                 \State \texttt{<Select random index>} 
                 \State \texttt{<Select random $x(0)$>} 
                 \State \texttt{<Select random trajectory length between 10 and 50}
                 \State \texttt{<Generate trajectory>}
                 \If {trajectory length $>$ 9}
                    \If{index $\leq$ 0.20}
                       \State \texttt{<Save test trajectory>}
                    \Else \If{index $\leq$ 0.35}
                       \State \texttt{<Save validation trajectory>}
                     \Else 
                       \State \texttt{<Save train trajectory>}
                     \EndIf
                \EndIf
                \Else 
                    \State \texttt{<pass>}
                \EndIf
                 \State \texttt{<i += 1>}
            \EndWhile
        \EndFor
    \EndFor
  \end{algorithmic}
\end{algorithm}

As described in \ref{alg:datasetcreation}, for every pair of ($\mu$, $\nu$) values, we generate 5 trajectories or random length (between 10 and 50), each one with a random initial condition chosen in each one of these 5 ranges [0.0, 0.2], (0.2, 0.4], (0.4, 0.6], (0.6, 0.8] and (0.8, 1.0]. We also have set some restrictions on the admissible values in a trajectory: no value in the trajectory can be smaller than $-1$ nor greater than $3$. If one of these conditions is satisfied, we remove the last point added to the trajectory and save the trajectory provided that it has a length greater than 9. Once a trajectory is generated, it is randomly assigned to the train, validation, or test data sets, resulting in a train-validation-test split with 1,110,266 trajectories in the train data set (65\%), 255,667 trajectories in the validation data set (65\%), and 341,374 in the test data set  (20\%).\medskip

Since the first layer in our neural network architecture is a convolutional layer, all trajectories have to be  padded (if needed) with 0s to reach the maximum length (50). In other words, 0s are inserted to the left of the trajectory to make all of them have a length equal to 50, see \cite[Ch. 5 \& Ch. 9]{goodfellow2016deep}. We also set a patience value equal to 20 to allow the training to end before the selected maximum number of epochs (200) if the training process does not improve the validation mean average error (MAE) after 20 consecutive epochs, thus saving training time. Finally, and for the sake of completeness, we have used a computer with 16 cores configured with 128 GB RAM and Nvidia RTX 3090 GPU with 22 GB RAM, running Ubuntu 20.10. The complete training process took less than 3 hours, running up to 24 epochs.\medskip

\section{Results}\label{sec:results}

\begin{figure}[ht]
 \centering
 \begin{subfigure}{.45\textwidth}
 \includegraphics[width=\linewidth]{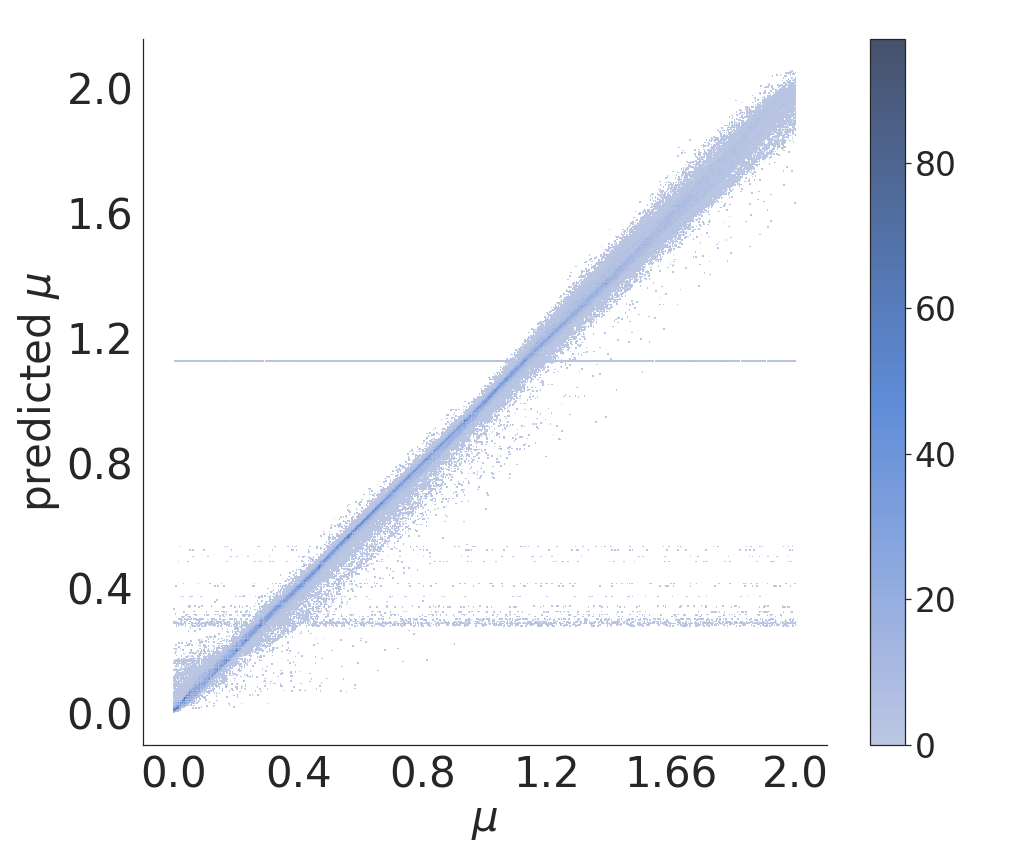}
 \end{subfigure}
 \begin{subfigure}{.45\textwidth}
 \includegraphics[width=\linewidth]{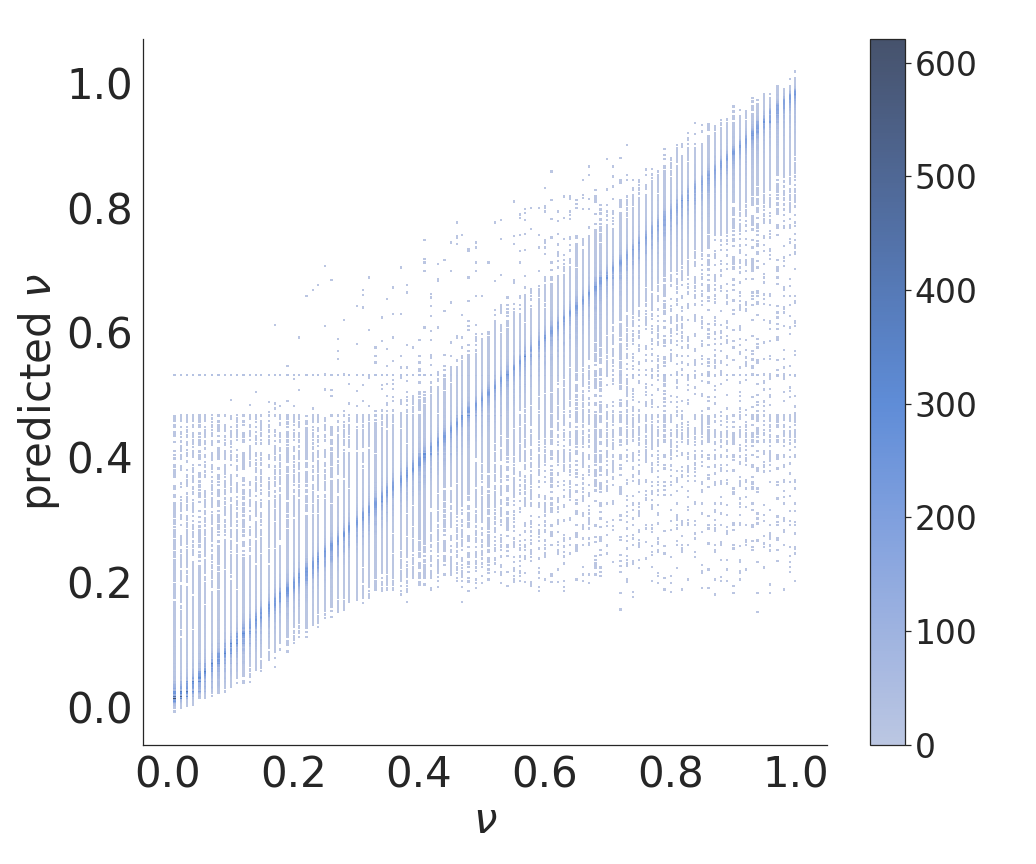}
 \end{subfigure} \\
 
 \begin{subfigure}{.45
 \textwidth}
 \includegraphics[width=\linewidth]{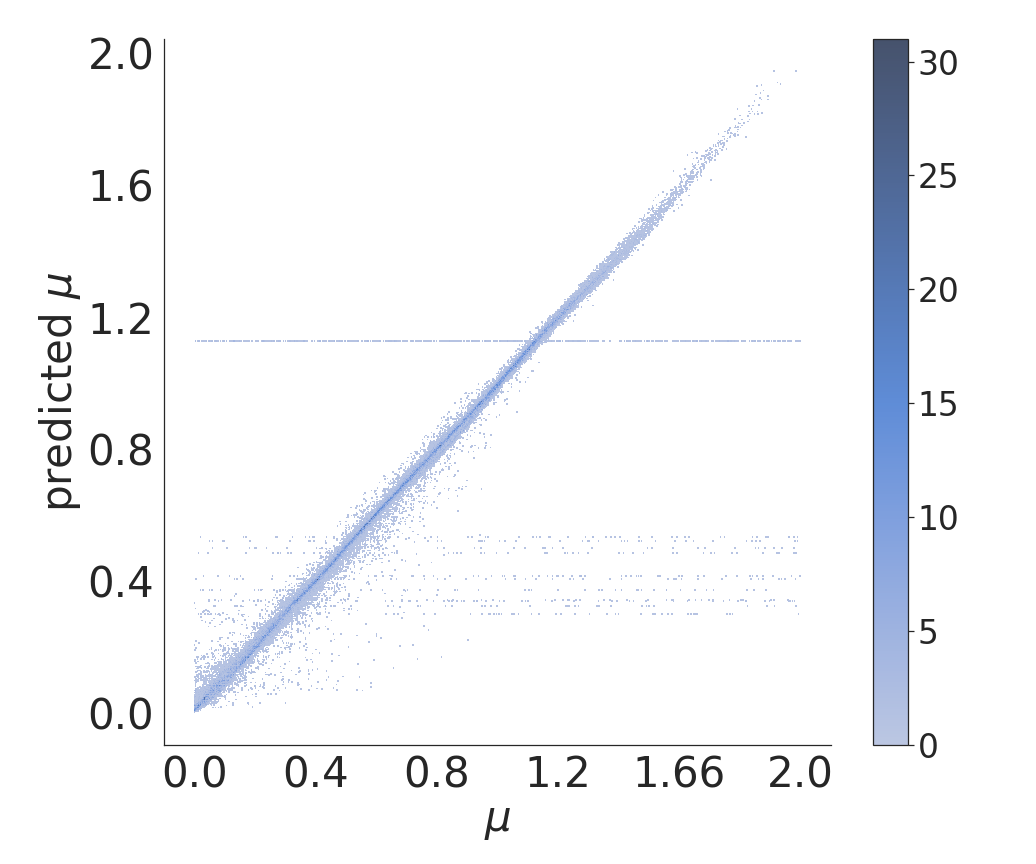}
 \end{subfigure}
 \begin{subfigure}{.45
 \textwidth}
 \includegraphics[width=\linewidth]{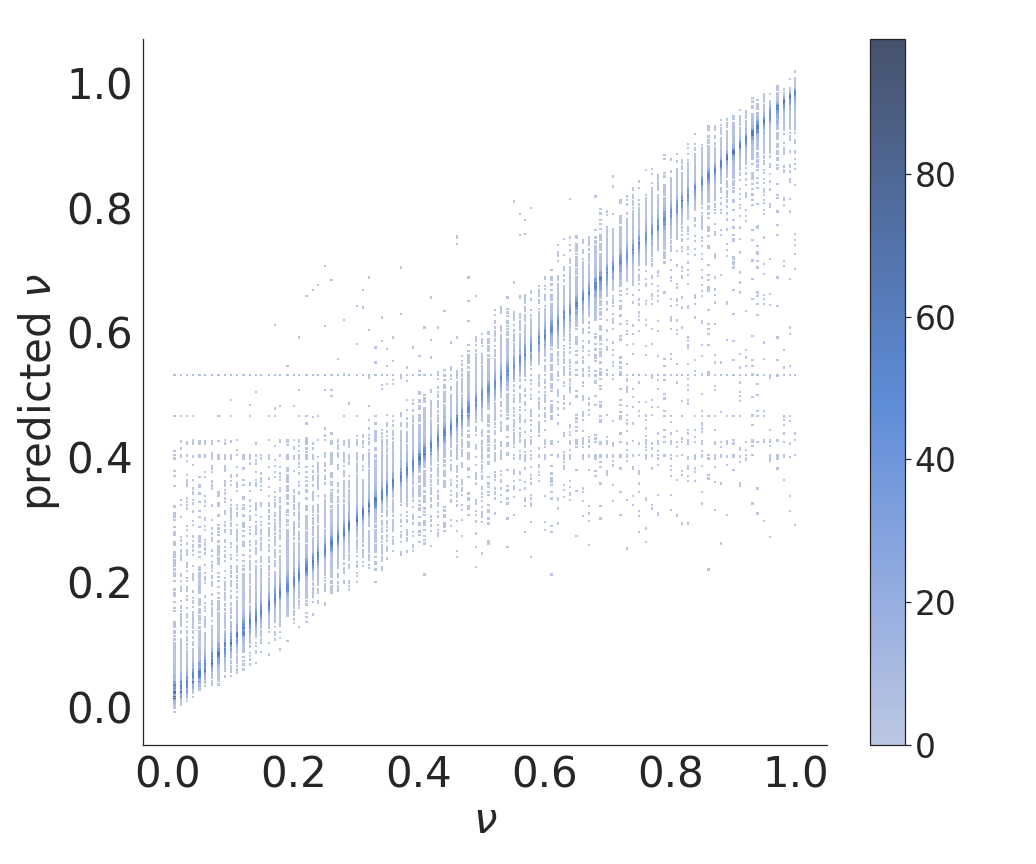}
 \end{subfigure}
 \caption{Truth vs predicted values of $\mu$ and $\nu$ in the validation data set (top) and for trajectories of length between 40 and 50 (bottom).}
 \label{fig:truth_vs_predicted2}
 \end{figure}

In Figure \ref{fig:truth_vs_predicted2}, we depict the MAE distribution for $\mu$ and $\nu$, comparing the ground truth (x-axis) and the predicted values (y-axis). The more values we have in the diagonal, the better the models can infer these values. The inference of the $\mu$ and $\nu$ parameters is very consistent along the entire range of possible values of each parameter. However, the MAE results for the inference of $\nu$ are slightly worse than the results for $\mu$. We also notice that the trajectory length affects the accuracy of the predictions. Since we are dealing with time series and recurrent neural networks, it is expected that long trajectories will return lower MAE values.\medskip

In Table \ref{table:mu_nu_MAE}, we show the MAE results for different lengths bins. MAEs both for $\mu$ and $\nu$ are in the same order of magnitude ($10^{-2}$). This is the order taken for the $\nu$ discretization. This also supports the use of a thinner discretization for $\mu$ ($10^{-3}$) respect to the one used for $\nu$.\medskip

MAEs for short trajectories are in fact considerably larger than for the rest of lengths
for instance, more than 20000 thousand for length equal to 10 but around 6000 for lengths between 40 and 50. As it was also noticed in \cite{conejero_garibo_lizama2022inferring} the MAE slightly increases for trajectories with lengths between $40$ and $50$. It is worth to mention, that as the length increases, the number of trajectories in each bin decreases due to the bounds set for stopping the generation of trajectories (no term can be between smaller than -1 nor greater than 3).\medskip

\begin{table}[t!]
    \centering
    \begin{tabular}{|c|c|c|c|c|}
    \hline 
    \bf Length & \bf MAE ($\mu$) & \bf MAE ($\nu$)\\ 
    \hline
    10-19 & 0.0284 & 0.0318\\
    20-29 & 0.0209 & 0.0229\\
    30-39 & 0.0188 & 0.0217\\
    40-50 & 0.0216 & 0.0227\\
     All  & 0.0233 & 0.0259\\
    \hline
    \end{tabular}
    \captionsetup{justification=centering}
    \caption{\label{font-table} $\mu$ and $\nu$ MAE's in the test data set.}
    \label{table:mu_nu_MAE}
\end{table}

\begin{figure}[t]
 \centering
 \begin{subfigure}{0.48\textwidth}
 \includegraphics[width=\linewidth]{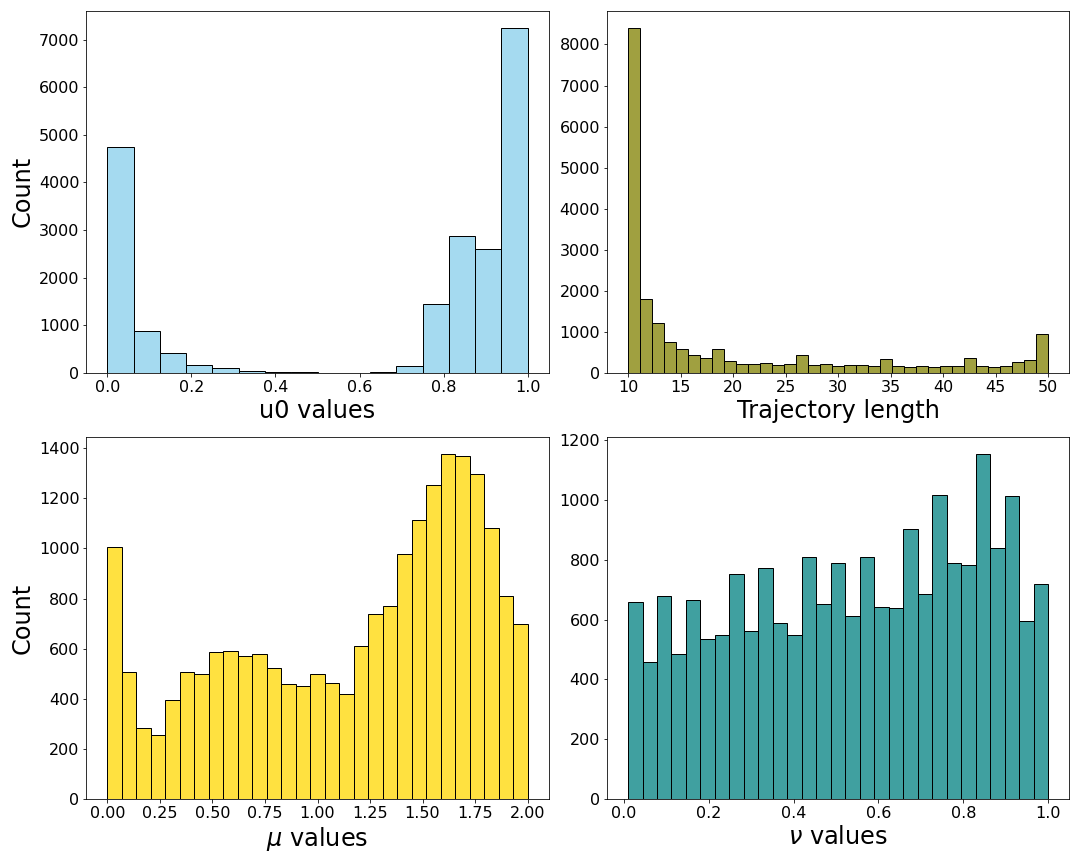}
 \end{subfigure}
\begin{subfigure}{0.48\textwidth}
 \includegraphics[width=\linewidth]{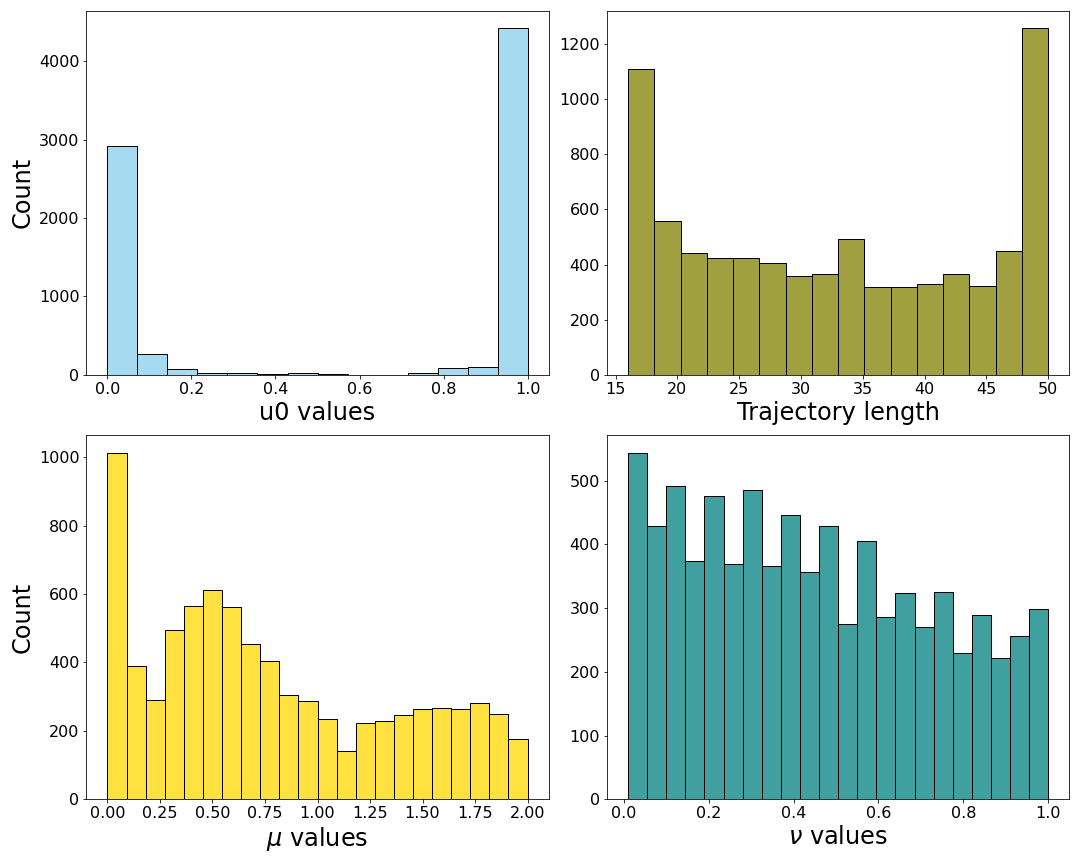}
 \end{subfigure}
  \caption{Histograms of $x(0)$, trajectory lengths, $\mu$, and $\nu$ values for trajectories with $\nu$ MAE higher than 0.05 (left) and for trajectories with $\nu$ MAE higher than 0.05 and length greater than 15 (right).}
   \label{fig:histograms_nu}
\end{figure}

In order to get deeper insights into the predictions obtained for the parameter $\nu$, we explore in more detail the ones that result in high values of MAEs. In Figure \ref{fig:histograms_nu}, we represent the frequency of trajectories with MAEs higher than 0.05 (left) in terms of the initial conditions $x(0)$, trajectory length, and the $\mu$ and $\nu$ parameters. We appreciate than in most cases, the main causes are short trajectory lengths, values of $\mu$ near to $0$ and an initial condition $x(0)$ close to 0.0 or to 1.0. On Figure \ref{fig:histograms_nu} (right), we show the same histograms but just for trajectories of length greater than 15. Here, we see that despite the trajectory length disappears, we still appreciate the same conclusions for the values of $x(0)$ and $\mu$. Besides, we also notice that whereas trajectories with values of $\mu$ between 0.4 and 1.1 return low MAEs for $\nu$, values of $\mu$ close to 0 result in high MAE values for $\nu$. 

\begin{figure}[t]
 \centering
 \begin{subfigure}{.40\textwidth}
 \includegraphics[width=\linewidth]{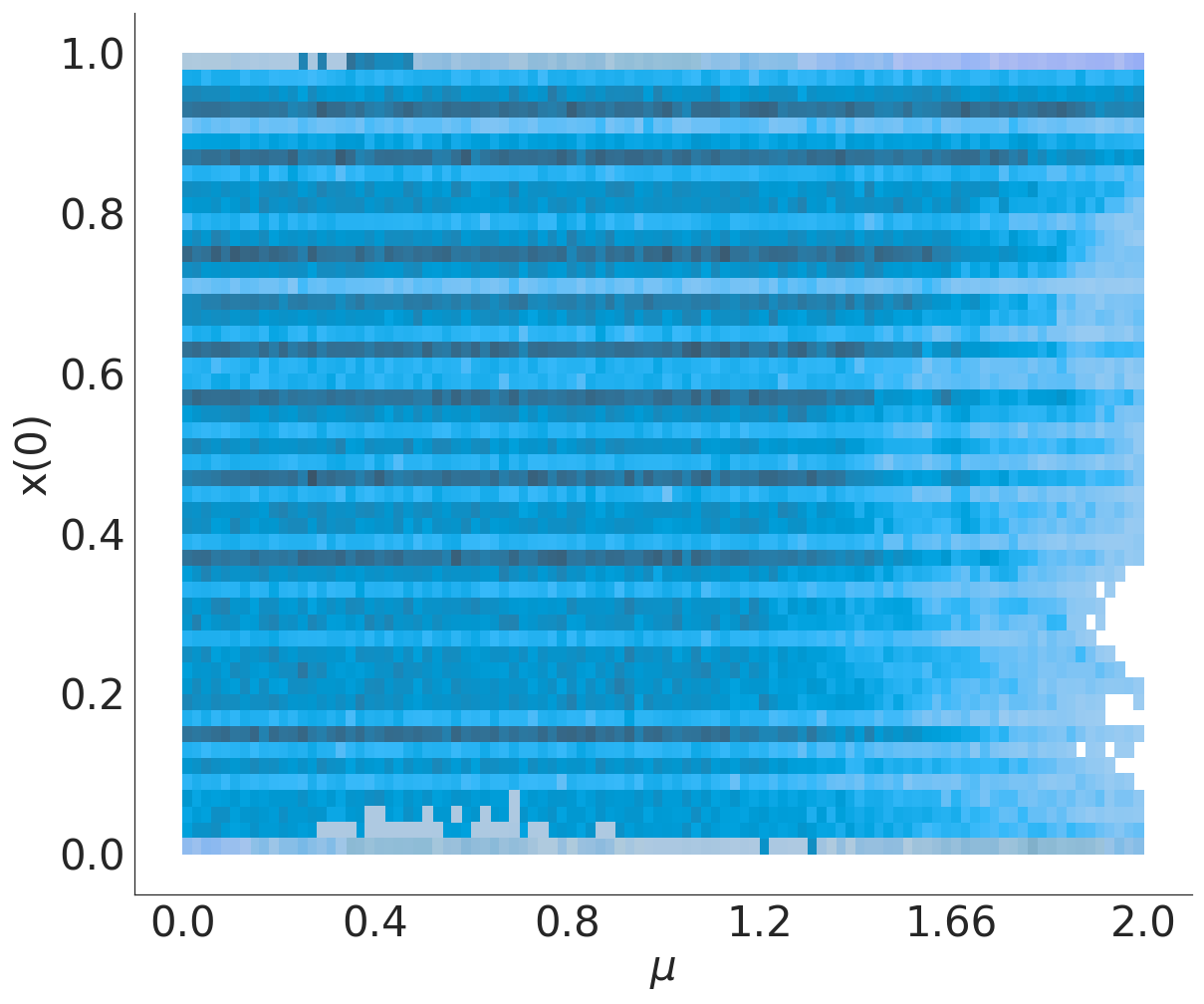}
 \end{subfigure}
 \begin{subfigure}{.05
 \textwidth}
 \includegraphics[width=\linewidth]{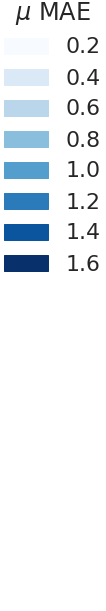}
 \end{subfigure}
 \begin{subfigure}{.40
 \textwidth}
 \includegraphics[width=\linewidth]{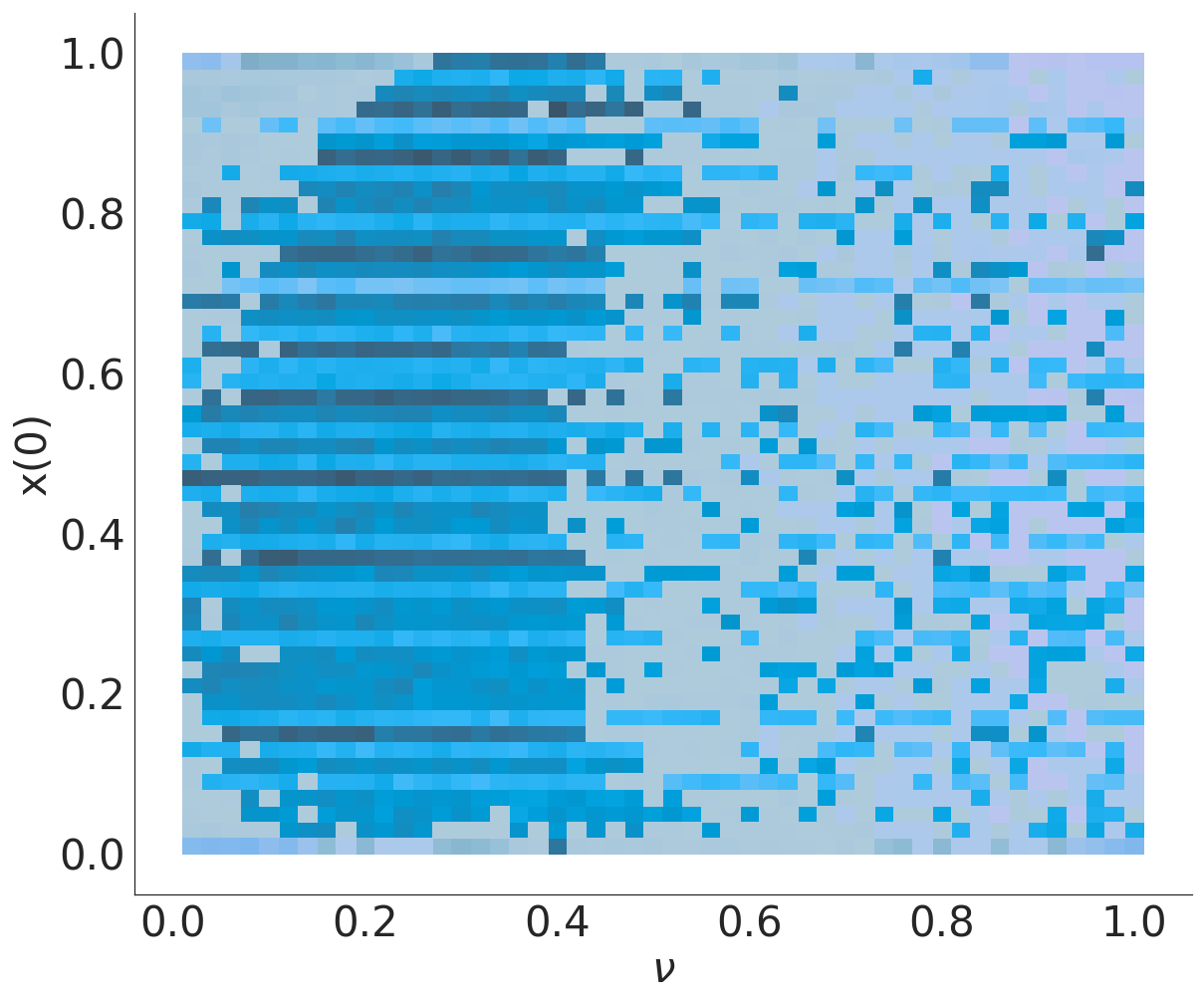}
 \end{subfigure}
 \begin{subfigure}{.05
 \textwidth}
 \includegraphics[width=\linewidth]{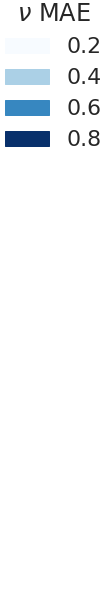}
 \end{subfigure}
\caption{$\mu$ (left) and $\nu$ (right) values vs $x(0)$ by $\mu$ and $\nu$ MAE respectively.}
 \label{fig:mu_nu_vs_u0}
 \end{figure}

We also confront the values of the initial conditions with the values of $\mu$ and $\nu$.
In Figure \ref{fig:mu_nu_vs_u0}, we represent the average MAE for pairs of initial conditions and values of the $\mu$ (left) and $\nu$ (right) parameters. The dark horizontal lines represent initial conditions that provide systematically higher values of MAE. On the one hand, for predictions of $\mu$, it seems that the initial condition is much more relevant for values of $\mu$ smaller than 1.66 than for greater ones. On the other hand, for the predictions of $\nu$, we can find initial conditions returning high values of MAE, mainly for values of $\nu$ smaller than 0.4.\medskip

\medskip

\begin{figure}[t]
	\centering
	\begin{subfigure}{.45\textwidth}
		\includegraphics[width=\linewidth]{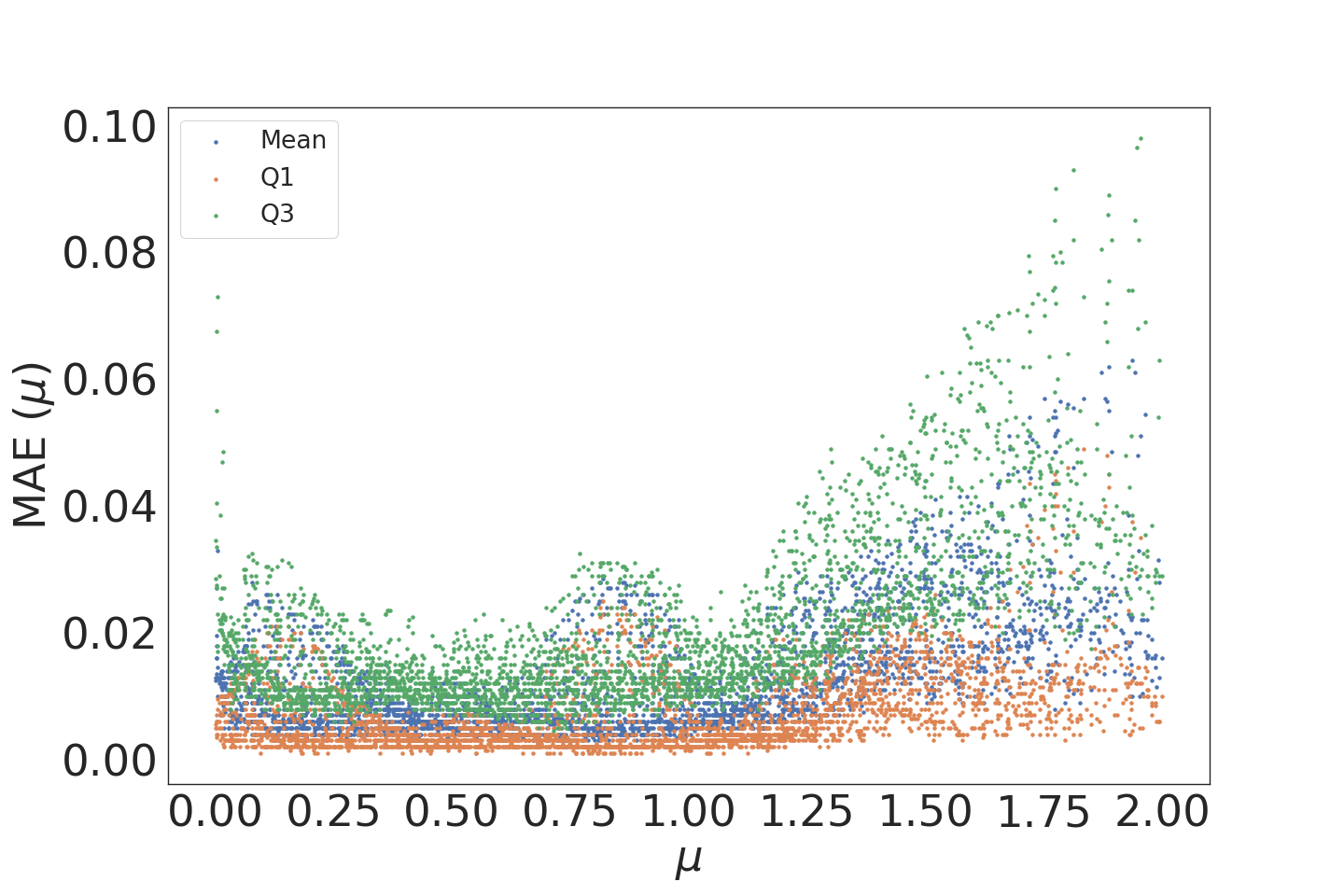}
	\end{subfigure}
	\begin{subfigure}{.45\textwidth}
		\includegraphics[width=\linewidth]{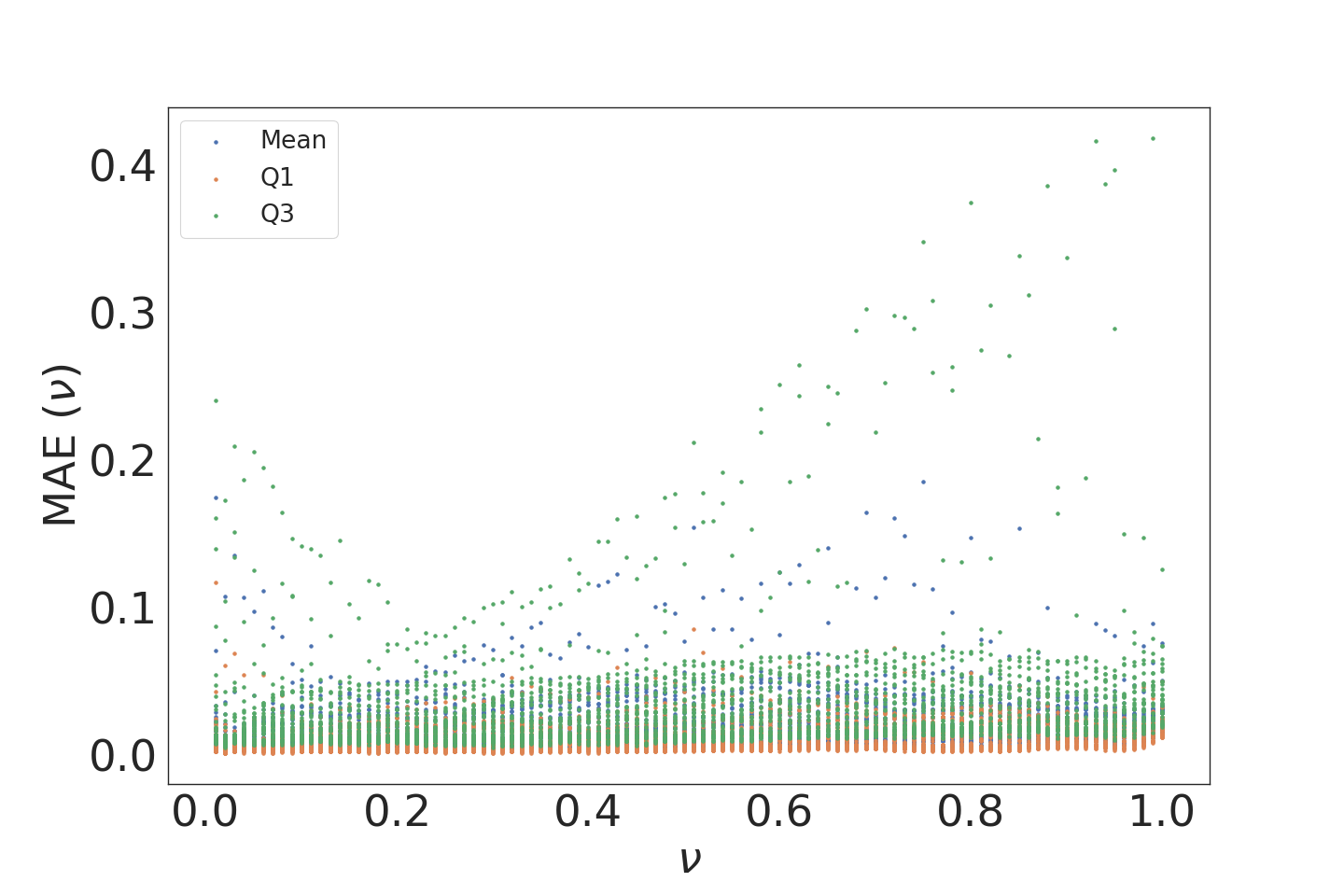}
	\end{subfigure} \\
	\caption{Quartiles Q1 (red), Q2 (blue), and Q3 (green) of the MAE distribution on the evaluation data for $\mu$ (left) and $\nu$ (right).}
	\label{fig:mae_q1_q2_q3_vs_mu_nu}
\end{figure}

We also study the MAE distribution for every $\mu$ and $\nu$ truth value. In Figure \ref{fig:mae_q1_q2_q3_vs_mu_nu}, we show the $Q_1$, $Q_2$, and $Q_3$ quartiles of these error distributions for each value of $\mu$ (left) and $\nu$ (right). On the one hand, it can be clearly seen that MAE error has greater variance for $\mu$ than for $\nu$, in particular for values of $\mu$ greater than 1.25, which is consistent with what we observed in Figure \ref{fig:mu_nu_vs_u0}. On the other hand, MAEs are considerably higher for $\nu$ rather than for $\mu$, as we have previously deduced from the Feigenbaum diagrams in Figures \ref{fig:feigenbaum_u0_03} and \ref{fig:feigenbaum_u0_03}, and from the comparison of predicted and truth values of $\mu$ and $\nu$ in Figure \ref{fig:truth_vs_predicted2}.\medskip

\begin{figure}[t]
	\centering
	\begin{subfigure}{.45\textwidth}
		\includegraphics[width=\linewidth]{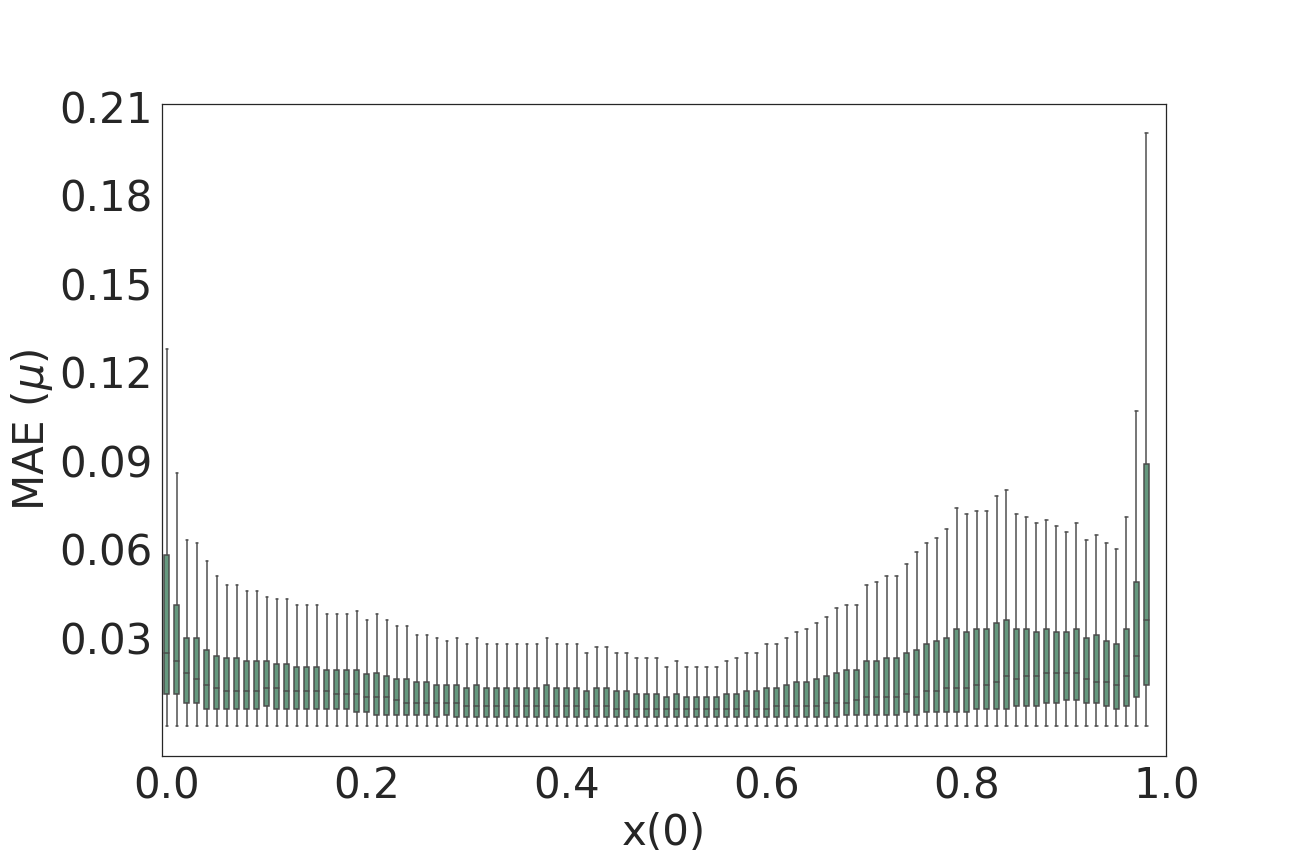}
	\end{subfigure}
	\begin{subfigure}{.45
			\textwidth}
		\includegraphics[width=\linewidth]{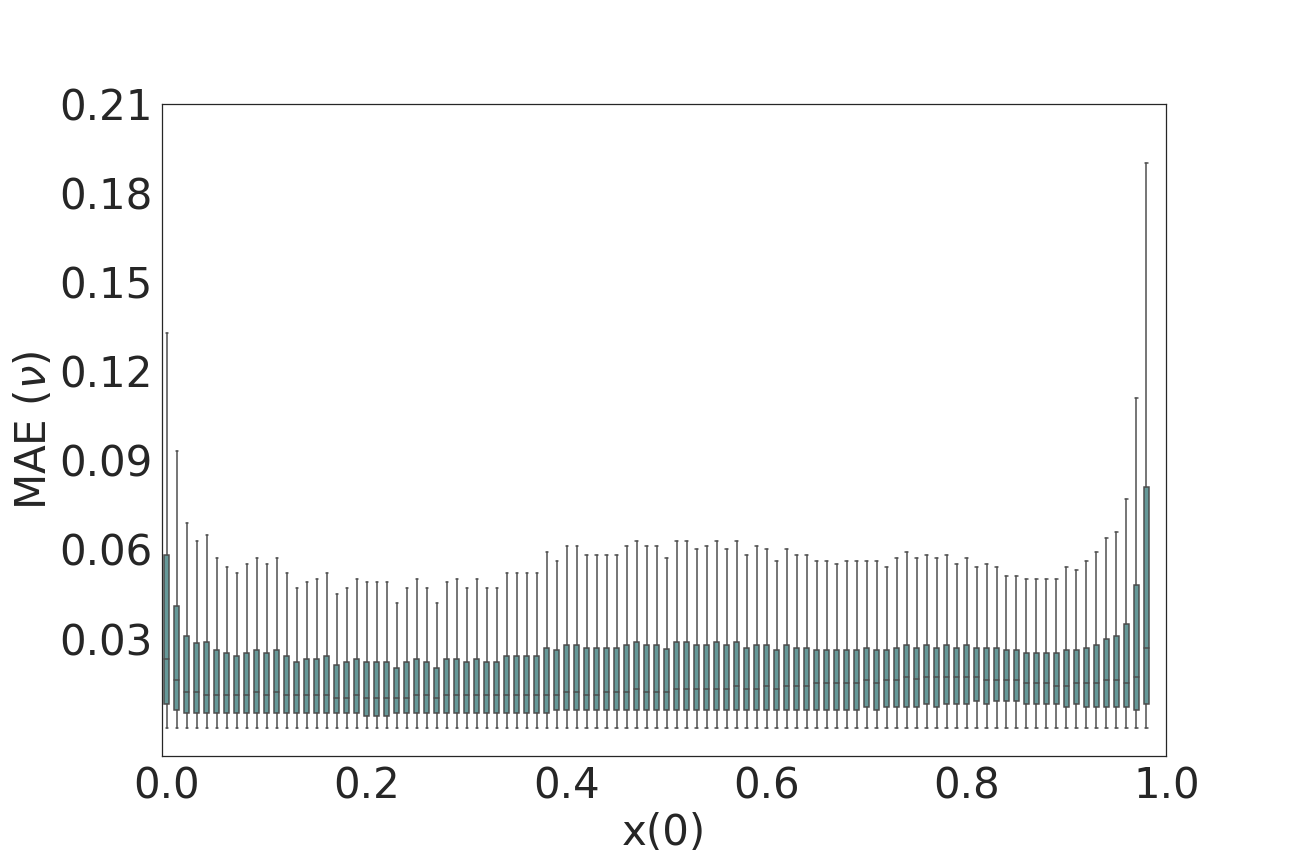}
	\end{subfigure}
	\caption{Box and whiskers plots for the MAE distribution obtained from predictions of $\mu$ (left) and $\nu$ (right) on the evaluation data set in terms of the initial condition $x(0)$. The boxes are painted in green and the whiskers on grey.} 
 \label{fig:boxplots_u0}
\end{figure}

Due to the fractional nature of equation \eqref{eq:wu_baleanu_delayed_formula}, the model has a memory component that is strongly dependent on the initial condition $x(0)$. In Figure \ref{fig:boxplots_u0}, we show boxplots of the MAE distribution for each initial condition $x(0)$.
The initial condition is more influential in predicting $\mu$ rather than in predicting $\nu$, since boxes (green) and whiskers (grey) are wider if we compare with the corresponding ones for $\nu$. Again, we see that the results are worse for initial conditions close to $0$ and to $1$, due in part to the form of the nonlinear logistic terms in \eqref{eq:wu_baleanu_delayed_formula}.\medskip

Finally, we wonder if machine learning methods are able to determine which trajectories are generated from \eqref{eq:wu_baleanu_formula} and which ones are given by \eqref{eq:wu_baleanu_delayed_formula}. In other words, we want to see if they can determine if a delayed is incorporated to the model or not.
To do so, we have used the data set described in this work jointly with the data set used in \cite{conejero_garibo_lizama2022inferring}. This second data set was also generated following the same procedure described in Algorithm \ref{alg:datasetcreation} but with the following specifications.

\begin{itemize}
	\item $\mu \in[2,3.2]$ discretized with a step size of $10^{-3}$,
	\item $\nu \in [0.01,1]$ discretized with a step size of $10^{-2}$,
	\item trajectory lengths between 10 and 50,
	\item and initial conditions $x(0)\in[0,1]$ randomly chosen with a resolution of $10^{-2}$.
\end{itemize}

\begin{figure}[t]
	\centering
	\begin{subfigure}{.45\textwidth}
		\includegraphics[width=\linewidth]{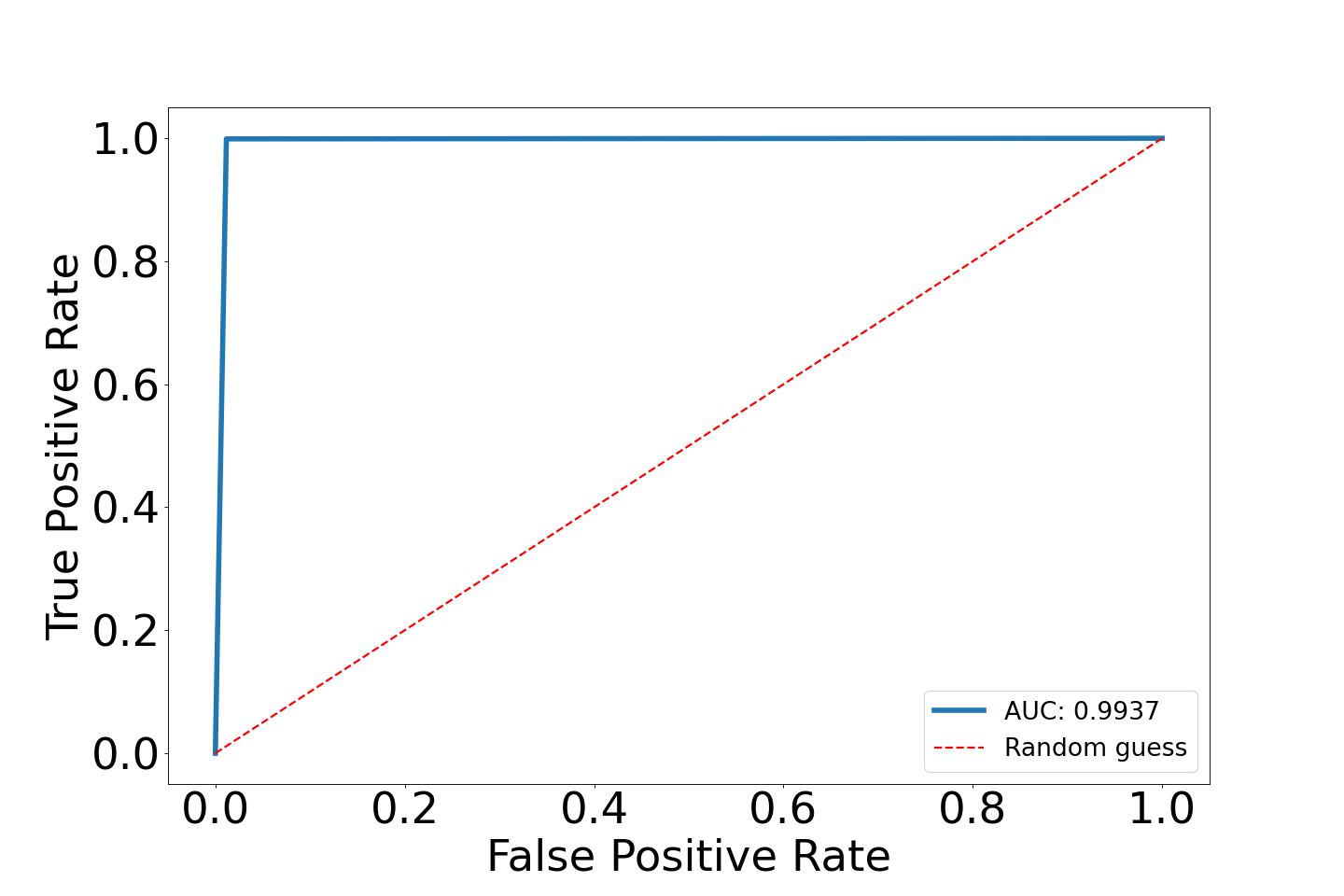}
	\end{subfigure}
	\caption{Receiver Operating Characteristic (ROC) curve on the test data set for the classification problem of fractional logistic trajectories with and without delay.}
	\label{fig:roc_auc}
\end{figure}

Since the range of values for $\mu$ is bigger in the present work than in the other, we have randomly sampled the train, validation and test data sets in order to obtain 3 perfectly balanced pairs of sets, consisting of 618,199 trajectories for each training data set, 142,800 for each validation data set, 190,334 for each the testing data set. For this classification task we have used the architecture  shown in Figure \ref{fig:architecture} but changing the activation function in the last dense layer to a sigmoid function. The training finishes after 43 epochs and 3 hours and 30 minutes because the 20 epochs of patience are reached. The accuracy in the validation data set is of a 99.38\%, which is very close to the area under the curve (AUC) value obtained for the test data set of 99.375\%, see Figure \ref{fig:roc_auc} where we show the receiver operating characteristic (ROC) curve and the AUC value.

\section{Conclusions}
\label{sec:conclusions}    

In this work, we have considered a delayed version of the fractional logistic equation introduced by Wu and Baleanu in \cite{wu_baleanu2015fractional}.
We have seen that a combination of convolutional and recurrent neural networks succeeds in training a model that given an new trajectory produced by this model is able to infer the $\mu$ and $\nu$ parameters of the model. We have seen that the MAE of the predictions falls almost in the order of magnitude of the discretization used. It would be interesting to check up to which point the results can be improved using other deep learning models, training with larger data sets, and with data sets of longer trajectories.\medskip

In our results, we have seen the importance of the initial condition, as the main contributor to the memory effect of the generating model. We have seen, that the initial conditions near $0$ and $1$ tend to accumulate predictions with higher errors in comparison with other values. However, we have detected some bins of initial conditions that produce higher errors than other neighbour bins, and this happens almost independently of the values of the parameters $\mu$ and $\nu$, see Figure \ref{fig:mu_nu_vs_u0}. We have also seen that the predictor of the parameters provides pretty accurate predictions for trajectory lengths greater than 15 and, in general, the accuracy of the predictions increases as long as the trajectory length increases.

Finally, we have seen that these methods almost perfectly detect if a trajectory comes from a particular model with or without delay.
It would be interesting to study whether this approach is able to determine the effect of the delayed term and of the fractional component in other fractional dynamical systems.
\medskip

\section*{Acknowledgements}
JAC acknowledges funding from grant PID2021-124618NB-C21 funded by MCIN/AEI/ 10.13039/501100011033 and by ``ERDF A way of making Europe'', by the ``European Union''.

\bibliographystyle{plain}
\bibliography{references}

\end{document}